\documentclass{article}
\usepackage{cite} 
\usepackage{tikz}
\usepackage{amsmath,amsthm, amsfonts}
\usepackage{color}
\usepackage[spanish,english]{babel}
\usepackage[latin1]{inputenc}
\usepackage{graphicx}
\usepackage{fancybox}
\usepackage{titlesec}
\usepackage{float}
\usepackage{enumerate}
\usepackage{multicol}
\usepackage{flushend}
\usepackage{pgfplots}
\usepackage[bookmarks]{hyperref}
\usepgfplotslibrary{polar}
\usepgflibrary{shapes.geometric}
\usetikzlibrary{calc}
\pgfplotsset{my style/.append style={axis x line=middle, axis y line=
middle, xlabel={$x$}, ylabel={$y$}, axis equal }}
\newtheorem{thm}{Theorem}[section]

\newtheorem{lem}[thm]{Lemma}
\newtheorem{prop}[thm]{Proposition}
\theoremstyle{definition}
\newtheorem{defn}[thm]{Definition}
\theoremstyle{remark}

\numberwithin{equation}{section}
\title{Asymptotic behavior of  global entropy solutions for nonstrictly
hyperbolic systems with linear   damping}
\author{Richard A. De la Cruz\\
Juan C. Juajibioy\\
Leonardo Rendón \\
    \small Bogotá}
\date{2014}

\begin{document}
\maketitle
\abstract{\noindent In this paper we investigate  the large
time  behavior  of the global weak entropy  solutions  to
the  symmetric  Keyftiz-Kranzer system with linear  damping. It  is
proved  that  as $t\to \infty$} the entropy solutions
tend  to zero  in the $L^p$ norm
\section{Introduction }\noindent
In this paper we  consider the Cauchy problem  to the symmetric system of Keyfitz-Kranzer type with linear  damping
\begin{equation}\label{intro-1}
 \begin{cases}
  u_t+(u\phi(r))_x+au=0,\\
  v_t+(v\phi(r))_x+bv=0.
 \end{cases}
 \end{equation}
with  initial  data 
\begin{equation}\label{intro-2}
 u(x,0)=u_0(x), \ v(x,0)=u_0(x), 
\end{equation}
This  system models of propagation  of forward longitudinal and
transverse waves of  elatic  string  wich moves in a plane, see \cite{kk1},\cite{cristes}. General source term for the system (\ref{intro-1})  was considered
in \cite{source1}. The damping  in the system (\ref{intro-1}) represents 
external forces  proportional to  velocity, and  this  term can  be produce 
lost of total energy of system. Consider the  scalar  case,by example

\begin{equation}\label{esc1}
 u_t+au_x+bu=0, \ u(x,0)=u_0(x).
\end{equation}
From the integral representation of (\ref{esc1}) it is  easy to find 
the following solution
\begin{equation}\label{esc2}
 u(x,t)=u_0(x-at)e^{-bt}.
\end{equation}
The term $bu$ produce  a dissipative  effect in the solutions,  i.e, the solutions
tends  to zero  when $t\to \infty$. We  are looking for  condition
under  wich the  terms $a$, $b$ have  a dissipative  efect in the solutions
of $\ref{intro-1}$.\\
Let $r(x,t)=\sqrt{u(x,t)^2+v(x,t)^2}$ be, we are  going  to show
the following main theorem.
\begin{thm}
If the initial  data $(u_0(x),v_0(x))\in L^{\infty}(\mathbb{R})\cap L^2(\mathbb{R})$
then  the Cauchy problem (\ref{intro-1}), ( \ref{intro-2}) has  a weak
entropy  solutions satisfaying 
\begin{equation}\label{in3}
 \|u\|_{L^{\infty}(\Omega)}+ \|u\|_{L^{\infty}(\Omega)}<M
\end{equation}
Moreover $r(u,v)$ converges to  zero in $L^{p}$ with  exponential
time decay, i.e.
\begin{equation}\label{in-4}
 \|r(x,t)|_{L^{p}(\mathbb{R})}\leq Ke^{-Mt} \|r(x,0)|_{L^{p}(\mathbb{R})}
\end{equation}
\end{thm}
\section{Preliminars}
We  start with some  preliminaries  about the general systems of conservation laws, see \cite{Bress1} chapter 5. Let
$f:\Omega\to \mathbb{R}^n$ be  a smooth  vector field. Consider  Cauchy problem for  the system
\begin{equation}\label{preli1}
 \begin{cases}
  u_t+f(u)_x=g(u),\\
  u(x,0)=u_{0}(x).
 \end{cases}
\end{equation}
When $g(u)=0$ the system (\ref{preli1}) is  called  homogeneous  system of  conservation laws, if $g(u)\neq 0$ the  system  (\ref{preli1}) is called
inhomogeneous system  or  balance  system  of  consevation laws. We shall work also  with  the parabolic perturbation to the system (\ref{preli1}), namely
\begin{equation}\label{preliparabol}
 \begin{cases}
  u_t+f(u)_x=\epsilon u_{xx}+g(u),\\
  u(x,0)=u_{0}(x).
 \end{cases}
\end{equation}
Denote  by $A(u)=Df(u)$ the Jacobian matrix  of partial  derivates of $f$.
\begin{defn}
 The system (\ref{preli1}) is strictly  hyperbolic  if for  every $u\in \Omega$, the matrix  $A(u)$ has  $n$ real
 distinct eigenvalues $\lambda_1(u)<\cdots <\lambda_n(u)$.
\end{defn}
Let  $r_i(u)$ the correspond  eigenvetor  to $\lambda_i(u)$, then
\begin{defn}
 We  say that  the i-th  characteristic field  is  genuinely  non-linear   if
 \begin{equation}\label{preli2}
  \nabla\lambda_i(u)\cdot r_i(0)\neq 0,
 \end{equation}
If instead 
 \begin{equation}\label{preli3}
  \nabla\lambda_i(u)\cdot r_i(0)= 0,
 \end{equation}
we say  that  the i-th characteristic field is  linearly  degenerate.
\end{defn}
\noindent For  the following  definitions  see \cite{smoller}, \cite{source2}
\begin{defn}
 A k-Riemann invariant is  a smooth  function $w_k:\mathbb{R}^n\to \mathbb{R}$, such  that 
 \begin{equation}\label{preli4}
  \nabla w_k(u)\cdot r_k(u)=0
 \end{equation}
\end{defn}
\begin{defn}
 A pair of function $\eta, q:\mathbb{R}^n \to \mathbb{R}$ is called  a entropy-entropy flux pair  if it  satisfies
 \begin{equation}\label{preli5}
  \nabla\eta (u)A(u)=\nabla q(u),
 \end{equation}
 if $\eta(u)$ is  a convex function then the pair $(\eta,q)$ is called  convex entropy-entropy flux pair.
\end{defn}
\begin{defn}
 A bounded measurable function $u(x,t)$ is  an entropy (or admisible) solution for  the Cauchy problem (\ref{preli1}), if it  satisfies  the following  inequality
 \begin{equation}\label{entropyineq}
  \eta(u)_t+q(u)_x+\nabla\eta(u)g(u)\leq 0. 
 \end{equation}
in the distributional sense, where $(\eta, q)$ is any  convex  entropy-entropy flux pair.
\end{defn}\noindent
We  consider the general system of  Keyftiz-Kranzer system
\begin{equation}\label{intro-1}
 \begin{cases}
  u_t+(u\phi(u,v))_x=0,\\
  v_t+(v\phi(u,u))_x=0,
 \end{cases}
 \end{equation}
to  get  some general observations  about this type of  systems.
Making $F(u,v)=(u\phi(u,v), v\phi(u,v)$  in (\ref{intro-1}), we have that the eigenvalues and  eigenvector  of the  Jacobian's matrix $Df$ are given by

\begin{align}
 \lambda_1(u,v)& =\phi(u,v)  & r_1&=(1,-\frac{\phi_u}{\phi_v}) \label{eq:1} \\
\lambda_2(u,v)& = \phi(u,v)+(u,v)\cdot\nabla \phi(u,v) &  r_2&=(1,\frac{v}{u}). \label{eq:2}
\end{align}
From (\ref{eq:1}),(\ref{eq:2}) we  have that $\nabla \phi\cdot r_1=0$,  and $\nabla Z(u,v)\cdot r_2=0$, where
$Z(u,v)=\frac{u}{v}$, then  the  Riemann invariants  are given by
\begin{align}
 &W(u,v)=\phi(r),\\
 &Z(u,v)=\frac{u}{v}.
\end{align}
\begin{lem}
 The system (\ref{intro-1}) is  always linear  degenerate in the  first  characteristic field. If
 \[
  (u,v)\nabla\phi(u,v)\neq 0,
 \]
then the system (\ref{intro-1}) is  strictly hyperbolic  
and non linear  degenerate in the  second characteristic field, moreover
\begin{equation}\label{r2condition}
 \nabla \lambda_2(u,v)\cdot r_(2)=\frac{2(u,v)\nabla\phi(u,v)+(u,v)H(\phi)(u,v)^{T}}{u}
\end{equation}
where $H$ represents the Hessian matrix.
\end{lem}

\begin{lem}
 Let $\eta(u,v) \in \mathbf{C^1}(\mathbb{R_+})$ a Lipschitz function in a neighborhod of the origin,  
 $q(u,v)=\psi(u,v)+\eta(u,v)\phi(u,v)$ be  a function, shuch that  $\psi$ satisfies
\begin{equation}\label{intro-3}
 \nabla \psi(u,v)=\left((u,v)\cdot \nabla\eta(u,v) -\eta(u,v)\right) \nabla \phi(u,v).
\end{equation}
Then the pair 
\begin{equation}\label{intro-4}
 (n(u,v), q(u,v))
\end{equation}
is a  entropy-entropy flux pair for the system (\ref{intro-1}). Moreover  if $\eta(u,v)$ is  a convex function, then the pair
(\ref{intro-4}) is  a convex  entropy-entropy  flux pair.
\end{lem}\noindent
\section{Global existence of weak entropy solutions and asymptotic
behavior}
We consider the parabolic  regularization of the system
(\ref{intro-1}), namely
\begin{equation}\label{g-1}
 \begin{cases}
  u_t+(u\phi(r))_x+au=\epsilon u_{xx},\\
  v_t+(v\phi(r))_x+bv=\epsilon v_{xx},
 \end{cases}
 \end{equation}
whit  initial data 
\begin{equation}\label{g-2}
 u^{\epsilon}(x,0)=u^{\epsilon}_0*j_{\epsilon}, \ v^{\epsilon}(x,0)=v^{\epsilon}_0*j_{\epsilon}, 
\end{equation}
where $j_{\epsilon}$ is a mollifier. 
In this  case $\phi(u,v)=\phi(r)$, with $r=\sqrt{u^2+v^2}$.
By (\ref{eq:1}) the  eigenvectors and eigenvalues  are given   by
\begin{align}
 \lambda_1(u,v)& =\phi(r)  & r_1&=(1,-\frac{u}{v}) \label{eq:11} \\
\lambda_2(u,v)& = \phi(r)+r\phi^{'}(r) &  r_2&=(1,\frac{v}{u}). \label{eq:22}
\end{align}
The following conditions will  be nesessaries in  our next  discution
\begin{enumerate}[$\text{C}_1$]
 \item $\lim_{r\to o}r\phi(r)=0$, $r\phi^{'}(r)\neq 0$
 \item $a>b$
\end{enumerate}
The condition $\text{C}_1$ garanties the strictly hyperbolicity to
the system (\ref{g-2}), while  condition $\text{C}_2$ ensure the existence
of  a positive invarian  region.
Now  we consider the  following subset of $\mathbb{R}$
\[
 \Sigma=\{(u,v): \phi(r)\leq C_0, 0<C_1\leq \frac{u}{v}\leq C_2\}.
\]
We affirm that $\Sigma$ is an invariant  region. Let $h(u,v)=(au,bv)$ be, if $(\overline{u},\overline{v})\in \gamma_1$ where
$\gamma_1$ is the level  curve of $Z=\phi(r)$ we have that
\[
 (\nabla W\cdot h)(\overline{u},\overline{v})=(a+b)r\phi^{'}(r)>0
\]
and if $\overline{u},\overline{v}\in \gamma_2$ where $\gamma_2$
is  the level curve of $Z$ we have that
\[
 (\nabla Z\cdot h)(\overline{u},\overline{v})=(a-b)\alpha_i>0
\]
with $i=1,2.$, then by the  Theorem 14.7 of \cite{smoller}, $\Sigma$
is an invariant  region  for the system (\ref{g-1}).Is easy  to verify
that $(au,bv)$ satisfies the condition $\text{H}_1\cdots  \text{H}_5$ in \cite{source1},
thus  we have the following  Lemma.
\begin{prop}
 If $(u_0,v_0)\in \Sigma$ and the $\text{C}$ conditions  holds,  then the Cauchy problem (\ref{g-1}),(\ref{g-2})
 has  a global  weak entropy solution.
\end{prop}
\noindent
Now  for the  global  behavior of  solutions, using  ideas of the author in \cite{panov2},
we  construct the following  entropy-entropy  flux pairs
\[
 n(r)=r^m, \ m\leq 2.
\]
From (\ref{intro-3}) we  have
\[
 q(r)=(m-1)\int_0^{r}s^m\phi^{'}(s)ds+r^m\phi(r),
\]
Integrating  by parts we have that
\[
 q(r)=(m-1)\int_0^{r}s^m\phi^{'}(s)ds+r^m\phi(r),
\]
integrating by parts  we have
\[
 q(r)=m\phi(r)-m(m-1)\int_0^r s^{m-1}\phi(s)ds.
\]
Let $M=\sup_{(u,v)\in [0,\|u|\_{L\infty}]\times[0,\|u|\_{L\infty}]}\{\phi(r)\}$, then
we have that
\begin{equation}\label{entro-1}
 |q(r)|\leq 2mM r^m.
\end{equation}
Multiplying in (\ref{intro-1}) by $\nabla \eta$ we have  that
\begin{equation}\label{entro-2}
 \eta(r)_t+q(r)_x\leq -3mMr^m
\end{equation}
Now  we choose $h(x)\in C^2(\mathbb{R})$ a function  such a
$|h^{'}(x)|\leq 1$, $|h^{''}(x)|\leq 1$ and $h(x)=|x| $ for $|x|\geq 1$
and  set $k(x)=e^{-h(x)}$, then $k^{'}(x)\leq k(x)$. Multiplying
by $k(x)$ in (\ref{entro-2}), and inegrating  over $x$ we have
\begin{equation}\label{entro-3}
 \frac{d}{dt}\int_{-\infty}^{\infty}\eta(r)g(x)\leq \int_{-\infty}^{\infty} q(r)k^{'}(x)+-3mM\int_{-\infty}^{\infty}r^mdx
\end{equation}
by the inequality (\ref{entro-1}) we have
\begin{equation}\label{entro-4}
 \frac{d}{dt}\int_{-\infty}^{\infty}\eta(r)k(x)dx\leq -mM\int_{-\infty}^{\infty} r^mk(x)dx.
\end{equation}
If $\psi(t)=\int_{-\infty}^{\infty}\eta(r)k(x)dx$  we have
\[
 \frac{d}{dt}\psi(t)\leq -mM\psi(t),
\]
by Gronwall's inequality  we have

\[
 \psi(t)\leq e^{-mMt}\psi(0).
\]
Thus we have
\begin{equation}\label{entro-5}
 (\int_{-\infty}^{\infty}r^m(t)k(x)dx)^{\frac{1}{m}}\leq e^{-Mt}(\int_{-\infty}^{\infty}r^m(0)k(x)dx)^{\frac{1}{m}}
\end{equation}
Passing to limit $m\to \infty$ in (\ref{entro-5}) we have the inequality (\ref{in3})

\section{Acknowledgments}
We would like to thanks to  professor Laurent Gosse by his suggestions and  review. To the  professor Juan Galvis   by his many valuable  observation, and to the professor Yun-Guang Lu by his suggestion this problem.
%-----------------------------------------------------------------
\bibliographystyle{amsplain}
\bibliography{state}
\end{document}